\documentclass[11pt]{article}
\usepackage{hyperref}
\title{Banach function algebras with dense invertible group}

\author{H.~G.~Dales and J.~F.~Feinstein}

\date{}

\usepackage{amsfonts}
\usepackage{latexsym}
\usepackage{theorem}

\theorembodyfont{\sl}
\newtheorem{theorem}{Theorem}[section]
\newtheorem{proposition}[theorem]{Proposition}

\def\C{{{\mathbb C}}}
\def\R{{{\mathbb R}}}
\def\T{{{\mathbb T}}}
\def\Q{{{\mathbb Q}}}
\def\N{{{\mathbb N}}}
\def\Z{{{\mathbb Z}}}
\def\F{{{\cal F}}}
\def\Re{{\rm Re}\,}
\def\Tint{{\rm int}\,}
\def\vz{{{\underline z}}}
\newcommand{\D}{\mathbb{D}}
\def\udisk{{\overline \D}}
\def\ri{{\rm i}}
\def\sr{{\rm sr}}
\def\tsr{{\rm tsr}}

\def\Inv{{\rm Inv}\,}
\newcommand{\fr}{\partial\,}

\newcommand{\ph}{\widehat}

\def\sskip{\vskip 0.3 cm\noindent}

\newcommand{\QED}{\hfill$\Box$}

\begin{document}

\maketitle

\begin{abstract}
In \cite{DawsonFein}, Dawson and the second author asked whether or 
not a Banach function algebra with dense 
invertible group can have a proper Shilov boundary. 
We give an example of a uniform algebra showing that this can happen, and 
investigate the properties of such algebras. 
\end{abstract}

\section{Introduction}

In \cite{Stolz}, Stolzenberg gave a counter-example to the conjecture
that, whenever a uniform 
algebra has proper Shilov boundary, 
its maximal ideal space must contain analytic structure
(see also \cite[Theorem 29.19]{Stout},
\cite{WermerStolz} and \cite{AlexanderWermer}). 
Cole gave an even more extreme example in his thesis \cite{C}, where
the Shilov boundary is proper and yet every Gleason part is trivial.

The requirement that the invertible group
be dense in the algebra is a stronger condition than the non-existence of
analytic structure, so a possible replacement conjecture is that
no uniform algebra with dense
invertible group can have a proper Shilov boundary. 

We modify the example of Stolzenberg to show that this 
conjecture is also false. Our example is also of the form $P(X)$ for
a compact set $X \subseteq \C^{\,2}$. In our example,
the Shilov boundary is
proper, and yet there is a dense set of functions in the algebra $P(X)$ whose
spectra have empty interiors. It is clear that this condition is
sufficient for the invertible group to be dense in the algebra. 

In fact, this sufficient condition is also
necessary for the invertible group to be dense, as is shown in \cite{FR}. 
Note, however, that it is not true that a function in the closure 
of the invertible group must have spectrum whose interior is empty.
For example, the invertibles are dense in  the uniform algebra 
$C[0,1]$, but there are functions in $C[0,1]$ whose range (and spectrum)
is the unit square.

Before constructing our example, we recall some standard notation and results.
In addition, we assume some standard knowledge of the theory of 
uniform algebras.
For more details see, for example, \cite{Dales},
\cite{Gamelin}, and \cite{Stout}.

Let $A$ be a unital algebra. Then $\Inv A$ denotes the invertible group of  $A$.
Now suppose that $A$ is a unital Banach algebra. Then we say that $A$ 
\textit{has dense invertibles}
if $\Inv A$ is dense in $A$.

In our terminology, a \textit{compact space} is a non-empty,
compact, Hausdorff topological space.
Let $X$ be a non-empty set. The uniform norm of a bounded,
complex-valued function $f$ on $X$ is denoted by $|f|_X$.
Now let $X$ be a compact space. 
The algebra of continuous, complex-valued functions on $X$
is denoted by $C(X)$. A \textit{uniform algebra on $X$} is a uniformly closed, unital
subalgebra of $C(X)$ that separates the points of $X$. Such algebras are, of course,
Banach algebras with respect to the uniform norm $|f|_X$.

Let $A$ be a uniform algebra on a compact space $X$. We denote by
$\Phi_A$ the maximal ideal
space of $A$; $A$ is \textit{natural} if $\Phi_A=X$. Recall that the 
{\it Shilov boundary} of $A$ is the minimum (non-empty) closed subset
$K$ of $X$ such that $|f|_K = |f|_X$ $(f \in A)$. The Shilov boundary of $A$
is denoted by $\Gamma_A$. For more details on the Shilov boundary, see
\cite[Chapter 7]{Stout} or \cite[Section 4.3]{Dales}.

The open unit disk in $\C$ is denoted by $\D$ and the boundary of a set 
$X \subseteq \C^{\,n}$ is denoted by $\fr X$.
For convenience, we shall
use the notation $\vz=(z,w)$ for a typical element of  $\C^{\,2}$.

Let $X$ be a compact subset of $\C^{\,n}$. 
The polynomially 
convex hull of $X$ is denoted by $\ph X$ and the
rationally convex hull of $X$ is denoted by $h_r(X)$.
The algebras
$P(X)$ and $R(X)$ are, respectively the uniform closures
in $C(X)$ of the set of restrictions to $X$ of
the polynomials and of the rational functions with poles off $X$.
It is standard that the maximal ideal spaces of $P(X)$ and $R(X)$ 
may be identified with
$\ph X$ and $h_r(X)$, respectively. With this identification,
$P(\ph X)$ and $R(h_r(X))$ are equal to the Gelfand transforms of the algebras
$P(X)$ and $R(X)$, respectively.

We are now ready to construct our example.

\begin{theorem}
There exists a compact set $Y \subseteq \fr \udisk ^{\,2}$ in $\C^{\,2}$
such that $(0,0) \in \ph{Y}$, and yet $P(Y)$ 
has dense invertibles. 
In particular, setting $X=\ph{Y}$,
the uniform algebra $A:=P(X)$ is natural on $X$ and 
$A$ has dense invertibles, but $\Gamma_A$ is a proper subset of
$X$.
\end{theorem}
{\bf Proof}
Let $\F$ be the set of all non-constant
polynomials $p$ in two variables with coefficients
in $\Q + \ri \Q$ such that $p(\udisk ^{\,2}) \subseteq \udisk$.
This is a countable 
family, and the family $\{p|Y:p \in \F\}$ 
is guaranteed to be dense in the unit ball of
$P(Y)$ whenever $Y$ is a compact subset of  $\udisk^{\,2}$.
We shall construct a compact set $Y \subseteq \fr \udisk ^{\,2}$ such that
$(0,0) \in \ph{Y}$ and such that, for all $p \in \F$, the spectrum of
$p|Y$ with respect to $P(Y)$ 
has empty interior. From this it quickly follows that $Y$ has the
desired properties.

Choose a countable, dense subset 
$\{\zeta_i : i=1,2,\dots\}$ of $\D$
which does not meet the countable set $\{p(0,0): p \in \F\}$.
Define sets $E_{i,p}$ for $i \in \N$ and
$p \in \F$ 
by 
\[
E_{i,p}=\{\vz \in \udisk^{\,2}: p(\vz)=\zeta_i\}\,.
\]
Each $E_{i,p}$ is compact, and there are only countably many such sets.
Enumerate those pairs $(i,p) \in \N \times \F$ for which $E_{i,p}$
is non-empty as $(i_j,p_j)_{j=1}^\infty$, and set
$K_j = E_{i_j,p_j}$ $(j \in \N)$.
For notational convenience, set $a_j=\zeta_{i_j}$, so that
$K_j=\{\vz \in \udisk^{\,2}: p_j(\vz)=a_j\}$.
Note that there may be repeats in the sequence $(a_j)$, and that the sets
$K_j$ need not be pairwise disjoint, but this will not matter.

Define polynomials $G_j$ for $j \in \N$ by
\[
G_j=\frac{p_j-a_j}{p_j(0,0)-a_j}\,.
\]
Then $G_j|K_j = 0$ and $G_j(0,0)=1$ for each $j \in \N$.

We now define inductively a sequence $(F_j)$ of entire functions 
which also satisfy the conditions that $F_j|K_j = 0$ and $F_j(0,0)=1$ 
for each $j \in \N$.

We begin with $F_1=G_1$.
Thus 
$F_1|K_1 = 0$
and 
$F_1(0,0)=1$.

Now suppose that $j\in \N$, and assume that 
$F_1,\dots,F_j$ have been defined so as to satisfy
$F_m|K_m = 0$ and $F_m(0,0)=1$ $(m=1,\dots,j)$.
Set 
\[
L_j = \{\vz\in\udisk^{\,2}:\Re F_j(\vz) \leq 1/2\}\,,
\]
and $H_{j+1}(\vz)=\exp(F_j(\vz)-1).$
We see that $|H_{j+1}|_{L_j} \leq \exp(-1/2) < 1$ and that 
$H_{j+1}(0,0)=1$.
Choose $N\in\N$ large enough that 
\[
|H_{j+1}^N G_{j+1}|_{L_j} < 1/4\,.
\]
Since $G_{j+1}|K_{j+1}=0$, we then also have
\[
|H_{j+1}^N G_{j+1}|_{L_j \cup K_{j+1}} < 1/4\,.
\]
Set $F_{j+1}=H_{j+1}^N G_{j+1}$.
Clearly $F_{j+1}|K_{j+1} = 0$ and $F_{j+1}(0,0)=1$.
The inductive construction may now proceed.

Note that this construction also produces a sequence $(L_j)$
of compact sets, where
 \[
 L_j =\{\vz\in\udisk^{\,2}:\Re F_j(\vz) \leq 1/2\} \subseteq \udisk^{\,2}\,.
 \]
We have $K_j \subseteq \Tint L_j$ $(j\in \N)$
(here and below the interior is taken in $\udisk^{\,2}$).
Moreover,
for each $j \in \N$, we have
\[
|F_{j+1}|_{L_j\cup K_{j+1}}< 1/4\,,
\]
and hence
$L_j \cup K_{j+1} \subseteq \Tint L_{j+1}$.

For each $j\in \N$, consider the variety 
$W_j=\{\vz \in \C^{\,2}: F_j(\vz)=1\}.$
Note that $(0,0)\in W_j$ and $W_j \cap L_j = \emptyset$.
Set $V_j=W_j \cap \udisk^{\,2}$. 
Since $W_j$ is the zero set of an entire function on $\C^{\,2}$ which vanishes
at $(0,0)$, it follows
(using, for example, 
\cite[Theorem C5, Chapter I]{GR})
that
$V_j \cap \fr \udisk ^{\,2} \neq \emptyset$.

Finally, set $M_j=\{\vz\in\udisk^{\,2}:\Re F_j(\vz) \geq 1/2\}$.
We now see that the following three facts hold.
\begin{enumerate}
\item[(a)]
The sets $\Tint L_j$ are nested increasing, and 
\[
\bigcup_{j=1}^\infty K_j \subseteq \bigcup_{j=1}^\infty \Tint L_j \,.
\]
\item[(b)]
Whenever $j \geq i$ we have $V_j \subseteq M_i$.
\item[(c)]
For each $j \in \N$, there is
a polynomial $h$ such that
\[
|h|_{M_j} < |h(\vz)|~~(\vz \in K_j)\,.
\]
\end{enumerate}
Here (a) and (b) are clear from the definitions and
properties above. 
For (c), note that
\[
|\exp(-F_j)|_{M_j} \leq \exp(-1/2) < 1 = \exp(-F_j(\vz))~~(\vz \in K_j)\,,
\]
and so we may take $h$ to be a suitable partial sum of the power series 
for $\exp(-F_j)$.

We now consider the space
of all non-empty,
closed subsets of $\udisk^{\,2}$, with the
Hausdorff metric.
(For more details on this metric,
see, for example, \cite[Appendix A.1]{Dales}.)
This  metric space is compact, 
so the sequence $(V_j)$ has a 
convergent subsequence, say
$V_{j_k} \to V \subseteq \udisk^{\,2}$.
We see that $(0,0) \in V$ because $(0,0)$ is in all of the sets $V_j$.
By property (b) above, $V \subseteq M_i$ for each
$i \in \N$.
Finally, as noted above, all of the sets $V_j$ have non-empty
intersection with $\fr \udisk ^{\,2}$, and so the
same is true of $V$.
Set $Y=V\cap \fr \udisk^{\,2}$, a non-empty, compact set.
We shall show that this set $Y$ has the desired properties.

We first show that $\ph Y = \ph V$ (and so, in particular, that
$(0,0) \in \ph Y$). For this it is enough to show that, for every 
polynomial $p$, we have $|p|_V = |p|_Y$.
Given a polynomial $p$, assume for contradiction that there is 
a $\vz_0 \in V \setminus \fr \udisk^{\,2}$ with
$|p(\vz_0)|>|p|_Y$.
Then there are disjoint open sets $U_1$ and $U_2$ in $\C^{\,2}$ with
$\vz_0 \in U_1 \subseteq \D^{\,2}$, with $Y \subseteq U_2$, and such that
$|p(\vz)|>|p|_{U_2}$ for all $\vz \in U_1$.
By the definition of $V$, there must be a $j\in\N$ such that 
$V_j \cap U_1 \neq \emptyset$ and $V_j \cap \fr \udisk^{\,2} \subseteq U_2$.
Thus $\max\{|p(\vz)|:\vz \in V_j\}$ is attained at a point of $\D^{\,2}$, and
not at any point of $V_j \cap \fr \udisk^{\,2}$. 
This contradicts the maximum principle (\cite[Theorem B16, Chapter III]{GR}) on $V_j$. 
Thus we must have
$\ph Y = \ph V$, as claimed.

Finally, let $p \in \F$. We wish to show that the spectrum of $p|Y$ in $P(Y)$
has empty interior, i.e., that the set $p(\ph{Y})=p(\ph{V})$ has empty interior.
Assume for contradiction that this is not the case. Then there must be
some $i\in \N$ with $\zeta_i \in p(\ph{V})$, and there is
some $\vz_0 \in \ph{V}$ with $p(\vz_0)=\zeta_i$.
In particular,
$\vz_0 \in E_{i,p}$, and so there is a $j\in\N$ with
$\zeta_i=a_j$, $p = p_j$, and $E_{i,p}=K_j$.
Thus $\vz_0 \in K_j \cap \ph{V}$.
Since 
$V \subseteq M_j$, it follows from (c), above, that
there is a polynomial $h$ with 
$|h|_V \leq |h|_{M_j} < |h(\vz_0)|$, and this contradicts 
the fact that $\vz_0 \in \ph{V}$.

We have proved that $Y$ has the desired properties.  Setting
$X=\ph{Y}$ and $A=P(X)$, the required properties of $A$ 
are now immediate from
the standard theory discussed earlier.
\QED
\sskip
{\bf Remarks}
Let $X$, $Y$ and $A$ be as in Theorem 1.1.
\begin{enumerate}
\item[(a)]
The coordinate projections $\pi_1:(z,w)\mapsto z$ and
$\pi_2:(z,w)\mapsto w$
are polynomials and are clearly in $\F$.
Thus the projections of $X$ on the two coordinate planes have empty 
interior, and in particular we must have $\Tint X = \emptyset$.
\item[(b)]
Denote by $\dim X$ the covering dimension of $X$. 
(See \cite{HW} for details of many equivalent
definitions of $\dim X$.)
Since $A$ is natural on $X$ and 
$A \neq C(X)$, we must have $\dim X > 0$.
By \cite[Theorem IV 3]{HW}, since $\Tint X=\emptyset$, we must have
$\dim X < 4$. Thus $\dim X \in \{1,2,3\}$. We do not know precisely
which values
of $\dim X$ are possible for our example.
\item[(c)]
Let $E$ be the projection of $X$ on the $z$-plane (i.e., $E=\pi_1(X)$).
Since $E$ has empty interior, the Shilov boundary of $R(E)$ is $E$.
By considering functions in $A$ of the form $f:(z,w)\mapsto r(z) p(w)$, 
where $r \in R(E)$ and $p$ is
a polynomial,
we see easily that
$X \cap (\udisk \times \T) \subseteq \Gamma_A$.
Applying the same argument to the second coordinate projection, we
obtain $X \cap \fr \udisk^{\,2} \subseteq \Gamma_A$, and hence we must have
$\Gamma_A = Y = X \cap \fr \udisk^{\,2}$.
\item[(d)]
If desired, it is easy to modify the construction above in order to ensure 
that $V$ (and hence $X$) is connected and is  
symmetric in the two coordinates.
\end{enumerate} 
\sskip

The following result shows that the example constructed above
has 
some further unusual properties.
This result is probably known, but we know of no explicit reference. 
For convenience, we use the term {\it clopen} to describe sets which are both 
open and closed.

\begin{theorem}
Let $A$ be a uniform algebra with compact maximal ideal space $\Phi_A$, 
and suppose that $A$ has dense invertibles. Then every component of every
{\rm (}non-empty{\rm )} zero set for $A$ meets $\Gamma_A$.
\end{theorem}
{\bf Proof}
Set $X=\Gamma_A$.
Assume, for contradiction, that $f \in A$ has a non-empty zero set $Z(f)$ 
which has a component 
$K \subseteq \Phi_A \setminus X$. 
Then $K$ is the intersection of the family of all relatively clopen subsets of $Z(f)$ which contain $K$, and so
there is a relatively clopen subset $E$ of $Z(f)$ with $K \subseteq E \subseteq \Phi_A \setminus X$.
We may then choose an open subset $U$ of $\Phi_A\setminus X$ with $E \subseteq U$ and such that
$Y:=\fr U \subseteq \Phi_A \setminus Z(f)$.
By the local maximum modulus principle (\cite[Theorem 9.8]{Stout}),
we must have $|g|_K \leq |g|_Y$ for all $g\in A$. 
Set $\delta = \inf \{ |f(x)|: x \in Y\} > 0$,
and then choose an invertible element $g \in A$ with $|g- f|_X<\delta/2$.  
Then  $|g^{-1}|_K > 2/\delta > |g^{-1}|_Y$, which contradicts the local maximum modulus principle.
\QED
\sskip

Let $X$ and $Y$ be as constructed in Theorem 1.1.
Let $\vz \in \ph X$, and suppose that $p$ is a polynomial with $p(\vz)=0$.
Then, since $P(X)$ has dense invertibles,
the zero set of $p$ must meet $Y$.
It follows that $\vz$ is in the rational hull $h_r(Y)$ of $Y$. 
Thus we have $h_r(Y)=\ph Y$ and $P(Y)=R(Y)$.

From the fact that $(0,0) \in \ph Y = h_r(Y)$ we quickly deduce 
that, for all $a \in \C$, there is some $(z,w) \in Y$ with $z=aw$. 
Set $U=\{(z,w) \in X: w \neq 0\}$. Then the image of $U$ under the rational
function $z/w$ is the whole of $\C$.
This suggests that the set $X$ is fairly large.

In fact, with a little more work, we can construct a fairly natural
continuous surjection from
a compact subset of $X$ onto a three-dimensional set.
In addition to $X$ and $Y$, let $V$ be as constructed in Theorem 1.1.
For $r$ and $s$ in $(0,1]$, set 
\[
B_{r,s}=\{(z,w) \in \C^{\,2}: |z|\leq r~\textrm{and}~|w|\leq s\},
\quad Y_{r,s}=V \cap \fr B_{r,s}
\]
and
\[
X_{r,s}=\ph{Y_{r,s}} \subseteq X\,.
\]
Using essentially the same arguments as 
above, one can show 
that 
$(0,0) \in X_{r,s}$, that the natural uniform algebra
$P(X_{r,s})$ has dense invertibles and that the
Shilov boundary of $P(X_{r,s})$ is $Y_{r,s}$.
Thus, as above,
for all $a \in \C$ there exists $(z,w) \in Y_{r,s}$ with $z=aw$.
Now set $E=\{(z,w)\in X: |w| \geq 1/4\}$ and consider the continuous map
$\eta: E \to \R \times \C$ given by
$(z,w) \mapsto (|w|,z/w)$.
Set $r=1$.
Let $s \in [1/4,1]$ and let $a \in \udisk$. Then, by above, there
exists $(z,w) \in Y_{r,s}$ with $z=aw$, and it is easy to see that we
must then have $|w|=s$.
From this it follows that 
$[1/4,1] \times \udisk \subseteq \eta(E).$

However it appears that we can not deduce anything about $\dim X$
directly from the fact that $\dim \eta(E) = 3$.
\sskip

We conclude this section with 
another example with some even stronger properties.

\begin{theorem}
There is a uniform algebra $A$ on a compact metric space
such that every point of $\Phi_A$
is a one-point Gleason part and such that
the invertible elements are dense in $A$, but
$\Gamma_A \neq \Phi_A$.
\end{theorem}

{\bf Proof} Taking the example above as our base algebra, 
we may form a system of
root extensions (as in \cite{C}) to obtain a new uniform algebra $A$ on a compact
metric space such that $\{f^2: f \in A\}$ is dense in $A$, and hence every point
of $\Phi_A$ is a one-point Gleason part.
Since the base algebra has proper Shilov boundary, the same is true for $A$ (\cite{C}).
Finally, since the base algebra has dense invertibles, so does $A$
(\cite[pp. 2837--2838]{DawsonFein}).\QED
\sskip
\section{Stable ranks}

We now discuss some related properties and questions concerning 
stable ranks. The reader may find the relevant definitions and 
standard theory
of the stable ranks discussed below 
in, for example, \cite{Badea}.
The connection with our work is that a unital Banach algebra has dense invertibles
if and only if it has topological stable rank equal to $1$.

There are many results known concerning
the Bass stable rank $\sr(A)$
and the topological stable rank $\tsr(A)$ of a unital Banach algebra $A$. 
(For more details on the history of this, and many more open questions, 
see, for example,
\cite{Badea} and \cite{CS}).
In particular, for all unital Banach algebras, we have
$\sr(A) \leq \tsr(A)$. (This inequality is strict for the disk algebra.)
If $A$ is a regular, commutative, unital
Banach algebra with
maximal ideal space $X$, then
\[
\sr(A) = \sr(C(X))=\tsr(C(X))=\left[\frac{\dim X}2\right]+1
\]
(where
$[t]$ denotes the greatest integer less than or equal to
$t$ and, as earlier, $\dim X$ is the covering dimension of $X$). 
Thus $C(X)$ has dense invertibles if and only if $\dim X \in \{0,1\}$.
In general, if $A$ is a commutative, unital Banach algebra 
with maximal ideal space $X$, we have $\sr(A)\leq\sr(C(X))$.
This inequality is strict in many cases, including the disk algebra.

If $A$ is a regular, commutative, unital Banach algebra
with maximal ideal space $X$, 
then the above facts give $\tsr(C(X))=\sr(C(X))=\sr(A)\leq\tsr(A)$
and so $1 \leq \tsr(C(X)) \leq \tsr(A)$.
In particular, if $A$ has dense invertibles
(so that $\tsr(A)=1$)
then we also have $\tsr(C(X))=1$, i.e.,
$C(X)$ has dense invertibles.

We shall strengthen this result slightly. Recall that 
a unital uniform algebra $A$ is {\it approximately regular} if every point of
$\Phi_A$ is an independent point. This is  equivalent to 
saying that every closed subset of $\Phi_A$ is $A$-convex.

Let $A$ be a unital Banach algebra. Then $\exp A = \{\exp a : a \in A \}$,
an open subset of $A$.

We shall  need the following standard characterization of the condition
$\dim X \leq 1$ for compact spaces $X$. This result is essentially
the case $n=1$ of
\cite[Theorem VI 4]{HW}.

\begin{proposition}
Let $X$ be a compact space. Then $\dim X \leq 1$ if and only if the 
following condition holds:
for every closed set $E \subseteq X$ and every $f \in \Inv C(E)$
there exists $\tilde f \in \Inv C(X)$ with $\tilde f | E = f$.
\QED
\end{proposition}

We are now ready to show that approximate 
regularity (rather than regularity) of $A$ is sufficient for
it to be true that whenever
$A$ has dense invertibles then $C(\Phi_A)$ has dense invertibles.

\begin{theorem}
Let $A$ be an approximately regular, unital  uniform algebra 
with maximal ideal space
$X$. Suppose that $A$ has dense invertibles. Then $C(X)$ has dense invertibles.
\end{theorem}

{\bf Proof}
Given that every closed subset of $X=\Phi_A$ is 
$A$-convex and that the invertibles are dense in $A$, we prove that the
covering dimension of $X$ is at most $1$ using the characterization given above: 
we show that for every closed $E \subseteq X$ and every
$f \in \Inv C(E)$, $f$ has a continuous extension in
$\Inv C(X)$. For this we use the Arens-Royden theorem
(see, for example, \cite[Chapter III, Theorem 7.23]{Gamelin}).

Take any such $E$ and $f$. 
Let $B$ be the closure of $A|E$ in $C(E)$. By assumption this is a natural 
uniform algebra on $E$.
By the Arens-Royden theorem, 
there exists $b\in \Inv B$ with $f \exp C(E) = b \exp C(E)$,
i.e., $f b^{-1} \in \exp(C(E))$.
Since  $A$ has dense invertibles,
we may choose an $a \in \Inv A$ such that
$a^{-1}|E$ is sufficiently close to $b^{-1}$ that 
$f \cdot (a^{-1}|E) \in \exp(C(E))$.
Choose $g \in C(E)$ such that $f a^{-1} = \exp g$ on $E$. Then $g$ has a 
continuous extension $\tilde g$ in $C(X)$. Set 
$\tilde f = a \exp \tilde g$. We have that $\tilde f \in \Inv C(X)$ and 
$\tilde f$ extends $f$,
as required.
\QED

\section{Open questions}

We conclude with some open questions. 

\begin{enumerate}
\item
Let $A$ and $X$ be as constructed in Theorem 1.1, and set
$M=\{f \in A: f(0,0)=0\}$. What can be said about $M^2$ and $\overline{M^2}$?
\item
Let $A$ be a uniform algebra with the property that
$\exp (A)$ is dense in $A$. Can the Shilov boundary of $A$ be proper?
Must $A$ be approximately regular? Must $A$ be $C(X)$?

It is conceivable (but unlikely) that $\exp A$ is dense in $A$ for the
example $A=P(X)$ of Theorem 1.1. 
This is true if and only if
$\exp A = \Inv A$, because $A$ has dense invertibles.
By \cite[Chapter III, Corollary 7.4]{Gamelin}, 
$\exp A = \Inv A$
if and only if $H^1(X,\Z)$, the first \u{C}ech cohomology group of
$X$ with integer coefficients, is trivial. Thus, for our example, 
$\exp A$ is dense in $A$ if and only if $X$ is simply coconnected in the sense
of \cite[Definition 29.24]{Stout}.
\item
Let $A$ be a uniform algebra with maximal ideal space $X$. Suppose that $C(X)$
has dense invertibles. Must $A$ have dense invertibles? Does it help to 
assume that $A$ is regular/approximately regular/has Shilov boundary $X$?
More generally (as asked in \cite{CS}), is it always true that
$\tsr(A)\leq\tsr(C(X))$?
\item
Let $A$ be a uniform algebra with maximal ideal space $X$.
Suppose that $A$ has dense invertibles. Must $C(X)$ have dense invertibles?
More generally (as asked in \cite{CS}), 
must we have $\tsr(C(X)) \leq \tsr(A)$?

Let $A$ and $X$ be as constructed in Theorem 1.1.
It is plausible that $\dim X >1$. If so, then we
have an answer to this question as
$C(X)$ would not have dense invertibles.

\item
Our final questions concern the existence of topological disks in the
maximal ideal space.
\begin{enumerate}
\item[(a)]
Let $A$ be a uniform algebra such that $\Gamma_A \neq \Phi_A$.
Does $\Phi_A \setminus \Gamma_A$ contain a homeomorphic copy
of $\udisk$?
\item[(b)]
Let $K$ be a compact subset of $\C^{\,n}$ such that $\ph K \neq K$.
Does $\ph K \setminus K$ contain a homeomorphic copy of $\udisk$?
\end{enumerate}
\end{enumerate}

\sskip
{\sf  
School of Mathematics

University of Leeds

Leeds LS2 9JT

England

email: garth@maths.leeds.ac.uk
\sskip

School of Mathematical Sciences

 University of Nottingham

 Nottingham NG7 2RD, England

 email: Joel.Feinstein@nottingham.ac.uk}
\sskip
2000 Mathematics Subject Classification: 46J10

\end{document}